\documentclass[a4paper, 12pt]{article}

\usepackage{amsmath}
\usepackage{amssymb}
\usepackage{amsfonts}
\usepackage{comment}
\usepackage{caption}
\usepackage{amsthm}
\usepackage[utf8]{inputenc}
\usepackage{graphicx}
\usepackage{hyperref}
\usepackage{url}
\usepackage{psfrag}
\usepackage{lmodern}
\usepackage{array}
\usepackage{wrapfig}
\usepackage{amsthm,amsfonts,amsmath,amscd,amssymb,verbatim}
\usepackage{setspace}

\usepackage[left=2.5cm, right=2cm, top=2cm, bottom=2cm]{geometry}

\theoremstyle{plain}
\newtheorem{theorem}{Theorem}[section]

\newtheorem{proposition}[theorem]{Proposition}

\theoremstyle{remark}

\newtheorem{remark}[theorem]{Remark}

\newtheorem{definition}[theorem]{Definition}
\newcommand{\dd}[2]{\frac{\partial {#1}}{\partial {#2}}}

\title{On an interaction model of general language change}
\author{A.~Fuchs, M.~Schwingenheuer, E.~Steinegger and T.~Voglmaier}
\date{\today}

\begin{document}

\maketitle

\abstract{In the following article, we construct an interaction model (a variant of the SIR-model) of general language change. In the context of language change it is desirable to deduce the long-term behaviour of the corresponding dynamical system (for example to decide if complete of reversible language change are going to happen). We analyse this dynamical system by first proving non-existence of periodic orbits and then invoking the Poincar\'{e}-Bendixson theorem to show convergence to critical points only. Non-existence of periodic orbits is established by contradiction in showing that the average position of a potential periodic orbit must coincide with a certain critical point $C$ which cannot be encircled by the flow of the dynamical system so that the average position would be pulled to that side. Thus the long-term behaviour of the model for any given initial constellation of speakers depends only on four interaction parameters and can be easily analysed by looking at the four critical points. Subsequent numerical analysis of real data on language change is used to justify the relevance of the constructed model for the practicing quantitative linguist. We show how data-fitting methods can be used to determine the four interaction parameters and predict from them the long-term behaviour of the system, i.e. if complete language change or reversible language change will take place.}
\thispagestyle{empty}
\tableofcontents
\thispagestyle{empty}

\section{Introduction and historical context}
Since its discovery, Piotrowski`s law \cite{P} has been tested and verified in various contexts of Quantitative Linguistics. The widespread use of this law, often called the "S-law``, has several reasons. Firstly, the similarity to the simplest mathematical model of the spread of contagious diseases, the SIR model (for an explanation see below), is evident (or more precisely, the Piotrowski-Altmann law can even be deduced from that analogy) and puts language change on a sound footing: the interaction of individual speakers. Secondly, the corresponding dynamical system is exactly soluble, thus there exists a mathematical function with two parameters (later three parameters) which can be fitted data-analytically to an empirically obtained dataset resulting in the possibility to immediately judge the quality of the modelling procedure. 
There has been some effort to apply the Piotrowski-Altmann law to situations like reversible language change, where for example a word initially gains popularity before later growing out of fashion (cf. the well-known paper by G.Altmann \cite{A}). Another situation where variants of the Piotrowski-Altmann law have been successfully applied is the two-stage change (Vulanovi\'{c} and Baayen \cite{VB}). Since Piotrowski's law \cite{P} is a priori either monotonically increasing or decreasing, it is principally not suitable for modelling such situations. In Altmann \cite{A} and Vulanovi\'{c} and Baayen \cite{VB}, the constant coefficient in the logistic differential equation is replaced by a linear or a higher-order polynomial in the time variable to account for reversible language change or the two-stage change. This is therefore of quantitative value but not of explanatory value for the underlying processes. As Altmann states in \cite{A} on pages 87-89, the next step should be a model of language change which comprises reversible language change in such a way that it can be traced back to the genuine interactions of speakers. This is the first step to connect language change to social, economic or historic reasoning. \\

The motivation of this paper lies in the fact that we wish to construct such an interaction model for language change (the so-called \emph{PLC model}) for speakers of some language, which contains Piotrowski-Altmann's law as a special case but which also contains other forms of language change like the reversible language change.

If we assume that the widespread use of the Piotrowski-Altmann law rests on the two criteria that it explains language change as an interaction process and that it is completely integrable, we meet the first of these two criteria. This is a starting point for the \emph{PLC model} to become a model of value to a practising quantitative linguist. Even though our model is not exactly soluble, it can be analysed with fairly straightforward mathematical methods, leading to a complete classification of the resulting scenarios of language change depending on four interaction parameters. The model is also accessible to numerical analysis, which we carry out doing in the last part of the paper to show how empirical data-fitting can be employed to obtain the interaction parameters in order to be able to make predictions about language change in special situations. We provide the used Code (in \emph{Python}) in the form of a \emph{jupyter-notebook} so that the interested reader has an interactive format to check the validity of the individual steps or to apply the described methods to own datasets.\newline
Crucial to showing the described classification of language change is the non-existence of periodic orbits of the underlying dynamical system, which we show in this paper. This is of mathematical interest in its own right. 
\newline
For the sake of completeness, we wish to point out that another approach, in spirit very similar to ours, is Wheeler \cite{LCCN}, who proposes similar techniques for a different model and with a different focus.

\section{The description of the PLC model}
\subsection{The motivational case: SIR model}
The PLC model (short for Progressive, Liberal, Conservative) is set up in analogy to the SIR model, which is often used to study the spread of contagious diseases like measles. Here $S$ stands for the number of "Susceptible Individuals", $I$ for the number of "Infectious Individuals" and $R$ for the number of "Resistant Individuals".
The SIR model is a system of coupled differential equations which mimick the interaction of the described groups. If for example a susceptible and an infected indiviual meet, there is a certain probability that the number of infected individuals increases by one and the number of susceptible individuals decreases by one (i.e. the susceptible individual has been infected). 
The model has gained some popularity in the Corona pandemic, cf. \cite{Pr}. We are interested in a slightly different model adopted to the current situation of language change, hence we are not going to delve into the SIR model, which has only been named for motivational reasons (the interested reader is referred to \cite{SIR}). Language change bears a strong resemblance: if two speakers meet (an ``infected'' one who already uses some new linguistic construct like a new word and a ``susceptible'' one, who is prone to learn the new word), there is a certain probability that the susceptible speaker will adopt the new word. We wish to point out that Piotrowski-Altmann's law can be deduced from a special case of the SIR-model (where no ``Resistant'' individuals are present) and results in the so-called ``logistic equation'', which arises in various areas of mathematical modelling.\\

\subsection{Assumptions of the PLC model}
In the following we discuss the assumptions made to obtain a meaningful model. Firstly we are, as in the case of the Piotrowski-Altmann law, taking no local variations like borders or other language boundaries into account. There are also models that do, but they are of a very different flavour, cf. \cite{VP}. We therefore assume a homogeneous speaker community of total size $N$ distributed evenly over some fictional country. We further assume that our country is divided in equally sized ``interaction'' spaces (imagine an equally spaced grid), meaning that if two individuals accidentally meet in such a space, there is a certain probability of interaction (expressed below by a Greek letter). Since we assume evenly distributed populations, the probability for an individual of group A to be in such an interaction space is proportional to the number $N(A)$ of individuals of population A. Assuming the independence of the spatial distribution of two groups A and B, the probability of an individual of A and one of B to meet in a certain interaction space is proportional to $N(A)\cdot N(B)$. 
We assume further that the total population of speakers has been exposed to a new linguistic construct like a new word which exists for a pre-existing one in the considered language (hence there is an old and a new version: \textit{ward, wurde} in German, cf. \cite{A}). 

Instead of the above-mentioned three groups (S,I,R) in the context of infectious diseases, we consider in the context of language change:

\begin{itemize}
 \item (P): progressive speakers, who only use the new feature (e.g. word) and try to spread its use.
 \item (L): liberal speakers, who are indifferent towards using the new feature (but do not speak it)
 \item (C): conservative speakers, who refrain from using the new feature and try to convince others to use the old feature.
\end{itemize}
As in the SIR model, $P(t)$ describes the number of individuals in the population of the progressive speakers at time $t$. Analogously for $L(t), C(t).$
From the assumption about the total size, it follows for all times: $P(t)+L(t)+C(t)=N.$

Now we describe the interactions:

\begin{enumerate}
 \item If a P and an L meet: There is a certain probability that L is converted to P (encoded in $\widetilde{\alpha}\geq 0$).
 \item If a C and an L meet: There is a certain probability that L is converted to C (encoded in $\widetilde{\gamma}\geq 0$).
 \item If a P and a C meet: There is a certain, individual probability that each is converted to L (encoded in $\widetilde{\beta},\widetilde{\delta}\geq 0$).
\end{enumerate}

\begin{remark}
 Please note that the Greek parameters introduced above are not actually the probabilities, since for the sake of simplicity additional aspects are incorporated into these parameters (e.g. number of interaction spaces and the like). Nevertheless they are \emph{proportional} to the interaction probabilities.   
\end{remark}

Thus the coupled system of differential equations of the PLC model reads:
\begin{align} 
P'&=&\widetilde{\alpha}  L P-\widetilde{\beta}  P C\\
L'&=&-\widetilde{\alpha} L P-\widetilde{\gamma}  L C +(\widetilde{\beta}+\widetilde{\delta})  P C\\
C'&=&\widetilde{\gamma}  L C-\widetilde{\delta}  P C \\
N&=&L+C+P
\end{align}

 A prime indicates differentiation with respect to time (the independent variable in the above system of coupled differential equations).
In the following, we are going to deduce properties of the described dynamical system and set them into perspective with regard to different scenarios of language change.

\section{Mathematical analysis of the model}
A priori, the \emph{PLC model} is three-dimensional, but since the total number of speakers is preserved, it can be reduced to a two-dimensional dynamical system. 
\subsection{Reduction to two dimensions}
Applying $N=P+L+C$ to eliminate the variable $L$, we obtain the system:

\begin{align}
P'&=&\widetilde{\alpha} P (N-P)-(\widetilde{\alpha}+\widetilde{\beta}) P C\\
C'&=&\widetilde{\gamma} C (N-C)-(\widetilde{\gamma}+\widetilde{\delta})  P C 
\end{align}

For convenience, it is best to normalise the involved quantities and express the system in terms of the lower-case letters and in the renamed parameters: 
\begin{align}
x=P/N;\ \ \ \ \ y=C/N;  \ \ \ \ \ \alpha=N \widetilde{\alpha};\\
\beta=N (\widetilde{\alpha}+\widetilde{\beta}); \ \ \ \ \ \gamma=N \widetilde{\gamma}; \ \ \ \ \ \delta=N(\widetilde{\gamma}+\widetilde{\delta}).
\end{align}

A short calculation shows that the \emph{PLC model} is equivalent to the following dynamical system:
\begin{align}
\dot{x}&=&\alpha  x (1-x)-\beta xy\\
\dot{y}&=&\gamma y(1-y)-\delta  x y
\end{align}

The dynamical variables $x,y$ describe the fraction of the total population of progressive and conservative speakers, respectively. Hence the range of $x,y$ is restricted to the unit interval $x,y \in [0,1]$ and $\alpha,\beta,\gamma,\delta \in \mathbb{R}_0^+$. Moreover, since the sum of the fractions cannot exceed 1: $x+y\leq 1$. Thus the dynamics of the system is restricted to the triangle: 
\begin{align}
\Delta=\{(x,y)\in \mathbb{R}^2 |0\leq x,y \ \wedge \ x+y\leq 1 \}.
\end{align}

In the sequel, we refer to the above dynamical system ($\Delta$ together with the equations) as the \emph{PC system}.

\subsubsection{Remark on the redefined parameters}
As can be seen directly from the definition $\alpha<\beta \  \Rightarrow \  \frac{\alpha}{\beta}<1$ as well as $\gamma<\delta \  \Rightarrow \  \frac{\gamma}{\delta}<1$ (if generically $ \widetilde{\alpha},\widetilde{\beta},\widetilde{\gamma},\widetilde{\delta}> 0$). Since this is crucial for the analysis in the following we generically deduce: $\alpha-\beta<0; \ \gamma-\delta<0$. Further, it follows straight away that $D:=\alpha\cdot\gamma-\beta\cdot \delta<0.$
For a discussion of the non-generic cases see section \ref{sing}.

\subsubsection{Remark on the relation to the Piotrowski-Altmann law}
What immediately springs to mind is the relation to the logistic growth. In each of the two equations, the first term precisely describes the logistic growth which in the second term is reduced proportionally to interactions.\newline
Thus, we can deduce instantaneously that the Piotrowski-Altmann law is comprised within the \emph{PC model}. To be more precise, if we start out with no conservative speakers at all, then their number remains zero for all times and we end up with:
\begin{align}
P'&=&\widetilde{\alpha} P (N-P)
\end{align}
which is the logistic differential equation or in relative terms:
\begin{align}\label{PAlaw}
\dot{x}&=&\alpha x(1-x)
\end{align}
For the sake of completeness, we include the solution to the differential equation \ref{PAlaw} (the usual form of the well-known Piotrowski-Altmann law):   
\begin{align}
x(t)=\frac{1}{1+e^{-\alpha t+b}}
\end{align}
To include incomplete language change, the law is generalised to:
\begin{align}
x(t)=\frac{c}{1+e^{-\alpha t+b}}
\end{align}
\subsubsection{Remark on possible generalizations and their implications}
It is not hard to show that the model can also be applied to certain situations where the interaction parameters are also negative. A necessary condition is given by:
\begin{align}
\delta+\beta \geq \alpha+\gamma
 \end{align}
 
This condition assures that the flow of the dynamical system points inwards of $\Delta$. 
The interpretation of the interactions is then as follows: A negative $\widetilde{\beta}$ for example can be interpreted as the possibility that progressive speakers can turn conservative speakers into progressive ones without having to change them to liberal speakers first. The actual process is more involved but can be seen in exactly that way. For example, if $\widetilde{\beta}$ is negative and $\widetilde{\delta}=-\widetilde{\beta}$, it is assured on the one hand that the number of liberal speakers will not change in time if progressive and conservative speakers meet. On the other hand, meeting of progressive and conservative speakers results in a decline in the number of conservative speakers and in an increase in the number of progressive speakers.   

Linguistically, this case may reflect a new loanword that does not seem really alien to a conservative speaker from the very beginning because its phonological and phonotactic structure is also typical of words of their own language, e.g. German \emph{cool, Trend, Flirt} or \emph{Trick}, taken in from English. (So in a way, the conservative speaker behaves progressively, possibly without intending to.)\newline
A very interesting issue arises if one allows more than one negative parameter. As it turns out in this case, given the right choice of parameters, the dynamical system even allows periodic orbits. This gives rise to a periodically repeating process! We plan to investigate this matter further in future research but refrain from doing so in the current paper.\newline  
A possible numerical example with one negative parameter (curve-fitting for the e-epithesis \ref{fig:5}) is given by:
$\alpha=0.1076; \ \beta=2.3732; \ \gamma=0.0377; \ \delta=-1.1806$. More details about the mathematical analysis in this case are given in section \ref{sing}. \newline

The methods described in the subsequent sections are also applicable to a much wider class of models describing the interactions of speakers. To give an inspiration of possible generalizations, we wish to point out that the conclusions of the paper carry over almost unaltered to models governed by systems of coupled differential equations of the form:
\begin{align}
\dot{x}=\alpha(x)(ax+by+C)\\
\dot{y}=\beta(y)(cx+dy+D)
\end{align}
where $\alpha,\beta$ are arbitrary positive, smooth functions and where $a,b,c,d,C,D$ are real numbers.

\subsection{The critical points}
A short calculation shows that the critical points of the PC system (where both $\dot{x}$, $\dot{y}$ vanish, i.e. where the dynamical system becomes stationary) lie at: 
\begin{align}
 C_0=(0,0); \ \ \ \ \   C_x=(1,0); \ \ \ \ \  C_y=(0,1); \\
 C=\left(\frac{\gamma\cdot(\alpha-\beta)}{D},\frac{\alpha\cdot(\gamma-\delta)}{D}\right).
\end{align}

It is not hard to see that the point $C$ always lies within $\Delta$. 
This corresponds to the following scenarios:
\begin{enumerate}
 \item At $C_0$: in our speaker community, there are only liberal speakers, but with no exposure to the new feature (so only the old feature is used).
 \item At $C_x$: in our speaker community, there are only progressive speakers, so everyone adopts the new feature (complete language change).
 \item At $C_y$: in our speaker community, there are only conservative speakers, so nobody adopts the new feature.
 \item At $C$: depending on the parameters, there is a fixed share of speakers using the new feature and the rest not using it (incomplete language change).
\end{enumerate}
\subsection{Outlook and putting into context}
We are going to show that in the long term, the scenarios described above are the only possible outcomes allowed by the \emph{PC model}. By this we mean that given the parameters $\alpha, \beta,\gamma,\delta$, any imaginable speaker constellation will for long times approach one of the described scenarios and we are able to tell what the outcome will be. We therefore expect great predictive power of our model when combined with suitable data-fitting procedures. To be more precise, we expect that given a dataset of how the usage of some new feature has developed over some period of time, after applying a data-fitting procedure to obtain estimates for the parameters, we can predict how the usage of the feature under consideration will terminate, i.e. if it becomes extinct, will be used by all speakers or only by a certain share of speakers. This gives a precise meaning to the classification of scenarios of language change claimed in the introduction. \newline

\subsection{Classification of long-term behaviour depending on generic parameters}
Dynamical systems, even in the supposedly easiest cases, show an intricate complexion of possible behaviour. The study of dynamical systems started with famous researchers like Newton, Lagrange and Poincar\'{e} trying to understand problems of celestial mechanics. Even two-dimensional dynamical systems can show chaotic behaviour and are in general hopelessly difficult to analyse (cf. Hilberts 16.th problem \cite{HI}). 
A recurring tool which often proves to be successful in the analysis of dynamical systems is the number of periodic orbits and their relation. A periodic orbit describes a specific configuration which returns to its initial position after evolving for some time $T$ (the period) according to the rules of the dynamical system. An infamous example in the solar system is given by the trajectory of the Earth around the Sun due to the gravitational law. Its period is obviously given by one year.
We are going to show in the following that the \emph{PC system} does not allow any periodic orbits. This will be the starting point for further analysis of the long-term behaviour of the \emph{PC system}. This section deals with the \emph{generic} case if all interaction parameters $\widetilde{\alpha},\widetilde{\beta},\widetilde{\gamma},\widetilde{\delta}$ are greater than zero. In section \ref{sing} we will discuss the model in a singular setting. In the following it is therefore understood that: $\widetilde{\alpha},\widetilde{\beta},\widetilde{\gamma},\widetilde{\delta}>0$.

\begin{theorem}\label{main}
The \emph{PC System} does not allow periodic orbits.
\end{theorem}

The idea of the proof of Theorem \ref{main} is as follows: we will argue by contradiction and assume that there exists such a periodic orbit $p(t)=(x(t),y(t))$ of period $T$. We proceed by showing that according to the rules of the \emph{PC system}, the average position 
\begin{align}
 \bar{x}=\frac{1}{T}\int_0^Tx(t)dt; \ \ \ \ \bar{y}=\frac{1}{T}\int_0^Ty(t)dt
\end{align}
of the periodic orbit coincides with the critical point $C$. The desired contradiction is then obtained by showing that a potentially existing periodic orbit cannot encircle the critical point but has to stay on one side of the critical point and would thus pull the average in $x,y$ away from $C$ towards that side.

\subsubsection{Fish-trapping in the \emph{PC system}}
In order to formalise the last step in the outlined sketch of the proof, we start by proving a proposition called the \emph{fish trap} in the following.
To this end, we firstly calculate the locus of vanishing derivative in the $x$-direction and the $y$-direction, respectively.
Setting $\dot{x}, \ \dot{y}$ to zero, a short calculation shows that 
\begin{align}
 g_x: \ y=\frac{\alpha}{\beta}(1-x); \ \ \ \ \ g_y: \ y=1-\frac{\delta}{\gamma}x,
\end{align}
are the lines of vanishing $\dot{x}$ and $\dot{y}$, respectively.

\begin{figure}[ht!]

\psfrag{w}{$\frac{\alpha}{\beta}<1, \ \frac{\gamma}{\delta}<1$}
\psfrag{x}{$g_x$}
\psfrag{y}{$g_y$}
\psfrag{c}{$C$}
\psfrag{1}{$I$}
\psfrag{2}{$II$}
\psfrag{3}{$III$}
\psfrag{4}{$IV$}
 \includegraphics[width=12cm]{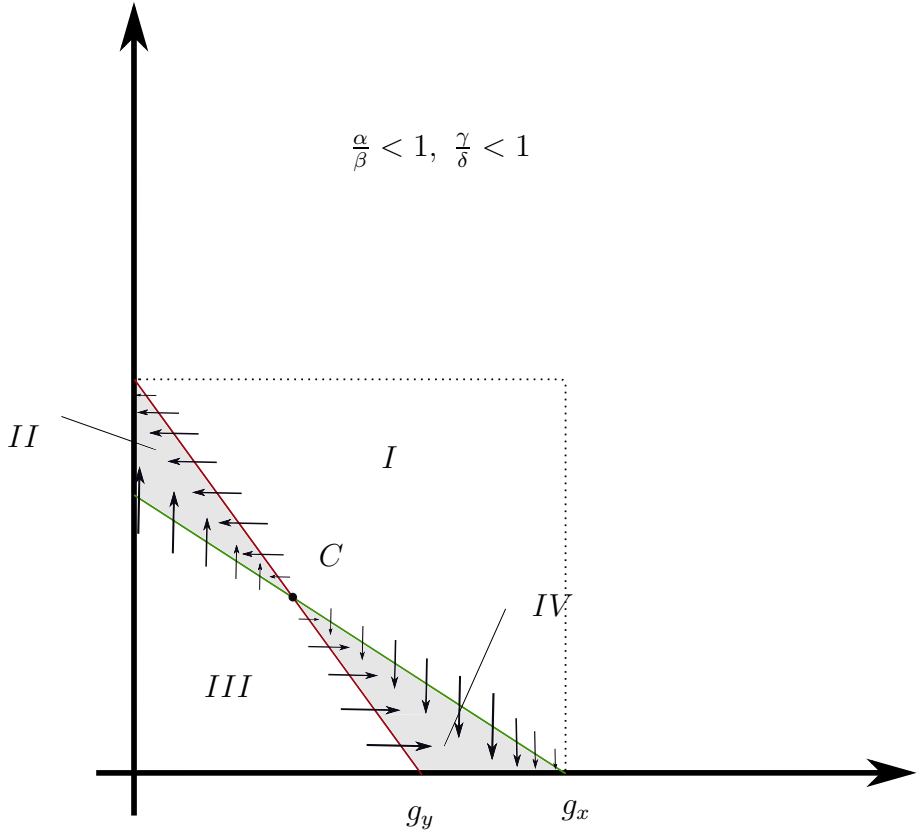}
 \caption{The fish-trap}
 \label{fig:1}
\end{figure}

Obviously, along the line $g_x$, the dynamical system has to flow in the $y$-direction (any flowing in the $x$-direction is prohibited by its very definition) and along the line $g_y$ the dynamical system flows in the $x$-direction. The intersection of both lines distinguishes the critical point $C$ at which all derivatives vanish. More can be said about the flow: above the line $g_x$, the $x$-derivative is negative, below the line $g_x$, the $x$-derivative is positive. Mutatis mutandis, the same is true for $g_y$. The area between the two lines together with the critical lines $g_x,g_y$ except point $C$ (cf. grey area in the Fig. \ref{fig:1}) make up the \emph{fish trap}, which traps any flow-line for all times (hence the nomenclature).

\begin{definition}
The sectors $II$ and $IV$ (Fig. \ref{fig:1}) together with the lines $g_x,g_y$ but without the critical point $C$ make up an area referred to as the \emph{fish trap}.
\end{definition}

\begin{proposition}[Fish trap]\label{fish-trap}
If a trajectory $\upsilon$ enters the \emph{fish trap}, it will stay there for all times. 
\end{proposition}

\begin{proof}
Without loss of generality, the trajectory $\upsilon(t)$ starts in sector $I$ and enters the fish trap on the boundary of sector $II$ at time $0$. Hence 
\begin{align}
 \upsilon(0) \in g_y; \ \ \ \ \dot{\upsilon}(0)=\left( 
\begin{array}{l}
-c \\ 
0 
                             \end{array}\right) \ \ \ \ c>0.
\end{align}

Now assume the opposite to the statement of the proposition, i.e. there exists a time $t$ such that $\upsilon(t)$ is outside the \emph{fish trap} either in sector $I$ or in sector $III$. By the Intermediate Value Theorem (IVT) \cite{B} there exists a time $\tau$ where $\upsilon$ leaves the \emph{fish trap} into the relevant sector, say without loss of generality, sector $I$. Thus $\upsilon(\tau) \in g_y$. To leave the \emph{fish trap}, the velocity vector $\dot{\upsilon}(\tau)=\left( 
\begin{array}{l}
\dot{x} \\ 
\dot{y}
\end{array}\right)$ has to point out of the \emph{fish trap} and thus has to satisfy:\\
\begin{minipage}{8cm}
\psfrag{c}{$C$}
\psfrag{1}{$I$}
\psfrag{2}{$II$}
\psfrag{3}{$III$}
\psfrag{4}{$IV$}
\psfrag{ft}{fish-trap}
\psfrag{g}{$\upsilon$}
\psfrag{v}{$\left( 
\begin{array}{l}
\delta \\ 
\gamma
\end{array}\right)$}
\psfrag{w}{$\left( 
\begin{array}{l}
\dot{x} \\ 
\dot{y} 
\end{array}\right)$}
 \includegraphics[width=5cm]{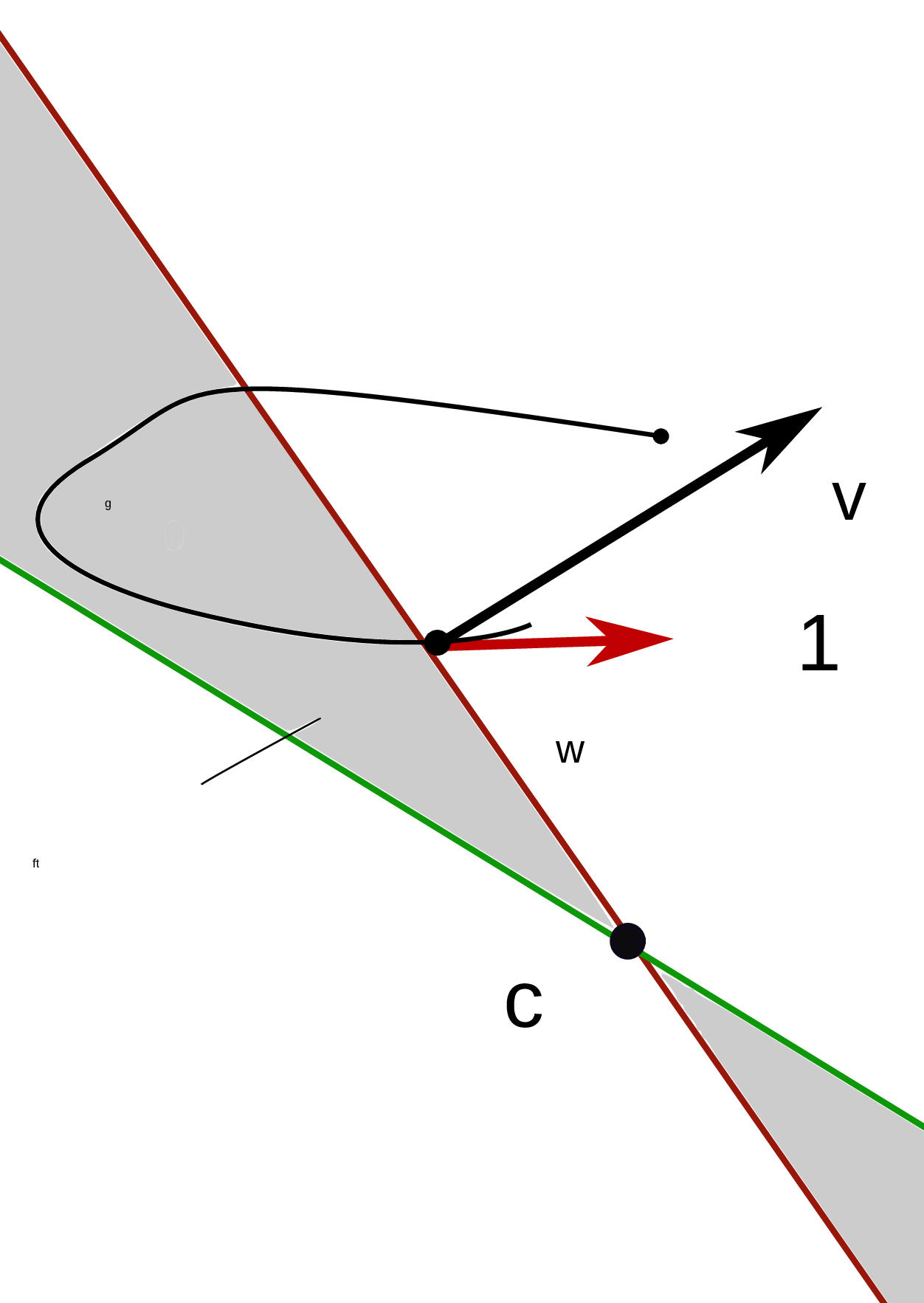}
 \label{fig:2}
\end{minipage}\begin{minipage}{3cm}
   $$\Rightarrow \left( 
\begin{array}{l}
\dot{x} \\ 
\dot{y}
\end{array}\right)\circ\left( 
\begin{array}{l}
\delta \\ 
\gamma
\end{array}\right)>0.$$
\end{minipage}
        
Hence $\delta \dot{x}+\gamma\dot{y}>0$. Now on the line $g_y$, it follows by definition that $\dot{y}=0$, hence $\delta \dot{x}>0$. With $\delta>0$, this yields $\dot{x}>0$. But this contradicts the fact that $\dot{x}$ is negative along that part of the line $g_y$ lying above $g_x$. Hence the assumption was wrong, and $\upsilon(t)$ cannot leave the \emph{fish trap} into sector $I$. The argument carries over verbatim to the case of entering into sector $III$ or with the trajectory starting in sector $III$. If it started already in sectors $II$ or $IV$, it could not leave by the same arguments. \end{proof}

\subsubsection{Non-existence of periodic orbits}
Now we can turn towards proving theorem \ref{main}.
\begin{proof}
Assume to the contrary that there exists a non-trivial periodic orbit of period $T$, which is notated by $p(t)=(x(t),y(t))$. By non-trivial we mean that it lies completely in the interior of $\Delta$ and that there exist times $t,t'$ such that $p(t)\neq p(t')$. In particular, no critical point lies on the trajectory of the periodic orbit, since it would have to stay there for all times. 
By the assumptions of non-triviality $0<x,y<1$, we can thus integrate 
\begin{align}
\int_0^T\frac{\dot{x}}{x}\ dt=\int_0^T(\alpha(1-x)-\beta y)dt.
\end{align}
Since the integral of $\frac{\dot{x}}{x}$ equals $\ln(x)$, by the periodicity, it follows straight away that the left hand side vanishes. Hence we find:
\begin{align}
 0=\alpha T -\alpha \int_0^Tx\ dt-\beta \int_0^Ty\ dt.
\end{align}
Dividing by $T$ yields:
\begin{align}
\alpha=\alpha  \bar{x}+\beta\bar{y},
\end{align}
where $\bar{x}=\frac{1}{T}\int_0^Tx\ dt$ denotes the mean value of $x$, analogously for $\bar{y}$.
The same can be done with $\frac{\dot{y}}{y}$, yielding:
\begin{align}
\gamma=\gamma  \bar{y}+\delta\bar{x}.
 \end{align}
Solving the two equations for $\bar{x},\bar{y}$, we can conclude that the point with coordinates $(\bar{x},\bar{y})$ coincides with the critical point $C$.
By non-triviality $P=p(0)$ is not any of the critical points, thus assume first that $P$ lies on or above the line $g_y$. 

If $p(t)$ stays above $g_y$ for all times (for example if it stays in sector $I$), then $y(t)>1-\frac{\delta}{\gamma}x(t)$ for all $t$. Integrating over the period $T$ and dividing by $T$ gives 
\begin{align}
 \gamma \bar{y}>\gamma-\delta \bar{x}
\end{align}
showing that the second equation obtained above for the mean values $\bar{x},\bar{y}$ is violated. Thus $p(t)$ has to cross $g_y$ eventually. By proposition \ref{fish-trap}, $p(t)$ is then confined to sector $II$ for all times coming. Since it cannot leave the fish trap anymore, being periodic, it already had to start there (it cannot reach the potential starting point anymore). But then it had been below $g_y$ already for all times, leading to an analogous contradiction as before. This shows that the assumption about the existence of a periodic orbit had been wrong from the beginning, proving the assertion of the theorem.  
\end{proof}
\subsubsection{Convergence of trajectories}
Due to H.~Poincar\'{e} and I.~Bendixson, there is a strong result about the behaviour of two-dimensional dynamical systems related to their periodic orbits. In the case we are concerned with, it basically asserts that any trajectory converges either towards a periodic orbit or to a critical point. Having previously excluded the existence of periodic orbits, using the Poincar\'{e}-Bendixson theorem, we can deduce strong conclusions for the \emph{PC system} about the long-term behaviour:

To state the theorem (and deduce our conclusions), we must first fix some notation.

\begin{definition}[Dynamical system]
Let $f \colon M \to \mathbb{R}^2$ be a smooth function from an open set $M \subset \mathbb{R}^2$ into $\mathbb{R}^2$. Then the differential equation 
\begin{align} 
\dot{\vec{x}}(t)=f(\vec{x}(t)),
\end{align}
defines a dynamical system on $M$.
\end{definition}

\begin{definition}[Trajectory of a dynamical system]
Given a dynamical system on the set $M$. A smooth curve $\upsilon \colon ]a,b[ \to M$ is called an integral of the dynamical system if its velocity vector $\dot{\upsilon}(t)$ coincides with $f(\upsilon(t))$ for all times. The set of points traced out in $M$ by $\upsilon$ is called the corresponding \emph{trajectory}.
\end{definition}
\begin{remark}
By abuse of notation, if there is no danger of confusion, we will sometimes use the terminology trajectory and integral interchangeably. 
\end{remark}

\begin{definition}[limit set]
Given a dynamical system on the set $M$.
A point $P\in M$ belongs to the \emph{limit set} of the point $Q\in M$ if there exists a sequence of times $t_j$ going to $\infty$ such that the integral $\upsilon$ of the dynamical system through $Q$ satisfies 
\begin{align}
\lim_{j \to \infty}\upsilon(t_j)=P.
\end{align}
We denote the \emph{limit set} (the set of all points sharing the property above) for $Q$ by 
\begin{align}
\omega(Q)\subset M.
\end{align}

\end{definition}

Now we can state:
\begin{theorem}[Poincar\'{e}-Bendixson]\label{PB}(Theorem 7.16 on page 223 in \cite{T})
Let a dynamical system on $M$ be given, fix a point $x \in M$ and suppose $\omega(x)\neq \emptyset$ is compact, connected and contains finitely many critical points. Then one of the following cases holds:
\begin{enumerate}
 \item $\omega(x)$ is a critical point
 \item $\omega(x)$ is a non-trivial periodic orbit
 \item $\omega(x)$ consists of finitely many critical points $\{x_j\}$ and non-closed trajectories connecting them.
\end{enumerate}
\end{theorem}

To deal with assumptions in theorem \ref{PB} we quote:

\begin{theorem}(Theorem 4.2 in \cite{V})\label{ass}
If a trajectory through $Q\in M$is confined to some bounded region, then the set $\omega(Q)$ is compact, connected and non-empty.
\end{theorem}

We need one more result for the proof:
\begin{proposition}\label{lim}
If $P \in \Delta$ is not a critical point, then there does not exist an integral $\upsilon \colon [0,\infty [ \to \Delta$ and a sequence $t_j \to \infty$ such that 
\begin{align}
\lim_{j \to \infty} \upsilon(t_j)=P.
\end{align}
\end{proposition}
\begin{proof}
Assume that there is such a point $P$. Since the argument carries over mutatis mutandis to the other sectors, we will assume without loss of generality that $P$ is in sector $I$. It is clear that $\upsilon$ has to stay in sector $I$ for all times, because once it would enter sectors $II$ or $IV$ it will never be able to return close to $P$ in sector $I$ due to proposition \ref{fish-trap}. But then $\upsilon$ will be above $g_y$ and $g_x$ for all times. By the defining property of $g_x,g_y$ this implies that for $\dot{\upsilon}=(\dot{x},\dot{y})$ we have: 
\begin{align} 
\dot{x}<0,\ \ \ \dot{y}<0.
\end{align} Hence $x(t_j)$ as well as $y(t_j)$ are monotonically decreasing sequences of real numbers bounded below by zero. By the Archimedean axiom and the assumptions made about $\upsilon$, they converge as follows: 

\begin{align} 
x(t_j)\to x(P); \ y(t_j)\to y(P).
\end{align}Note that it must hold 
\begin{align} 
x(t_j)\geq x(P); \ \ y(t_j)\geq y(P) \ \ \forall j.
\end{align}Due to the monotonicity it also follows for 
\begin{align} 
t_j<t<t_i \ \Rightarrow \ x(P)\leq x(t_i)\leq x(t)\leq x(t_j).
\end{align} 
Hence for $t>t_j$ it follows that 
\begin{align} 
||\upsilon(t)-P||^2= |x(t)-x(P)|^2+|y(t)-y(P)|^2\leq\\
|x(t_j)-x(P)|^2+|y(t_j)-y(P)|^2, \\
\Rightarrow \ \  ||\upsilon(t)-P||\leq ||\upsilon(t_j)-P||.
\end{align}
Now we are going to evaluate the mean values as in the proof of theorem \ref{main} and deduce the desired contradiction.
Given any $\epsilon>0$, choose a $\nu>0$ such that for $x$ within distance $\nu$ of $x(P)$ and $y$ within distance $\nu$ of $y(P)$ we have:
\begin{align}
 |ln(x)-ln(x(P))|<\epsilon; \ \ |ln(x)-ln(x(P))|<\epsilon.
\end{align}

Such a $\nu$ exists by the continuity of $ln$ in $x(P),y(P)$. 

Now choose $j$ so big such that 
\begin{align}
 ||\upsilon(t_j)-P||<\nu,
\end{align}
hence we find for all $t>t_j$:
\begin{align} 
|x(t)-x(P)|<\nu \Rightarrow \ |ln(x(t))-ln(x(P))|<\epsilon \\
|y(t)-y(P)|<\nu \Rightarrow \ |ln(y(t))-ln(y(P))|<\epsilon.
\end{align}

As before consider the integrals 
\begin{align} 
\left| \frac{1}{t-t_j}\int_{t_j}^t\frac{\dot{x}}{x}dt\right|=\frac{1}{t-t_j}\cdot\left| ln (x(t))-ln(x(t_j))\right|\leq \frac{\epsilon}{t-t_j}.
\end{align}

On the other hand, we have from the \emph{PC system}:
\begin{align}
\frac{\dot{x}}{x}=\alpha (1-x)-\beta y.
\end{align}
Integrating as before shows 
\begin{equation}\label{eq:inq1}
 \left|\alpha-\alpha \frac{1}{t-t_j}\int_{t_j}^tx dt-\beta \frac{1}{t-t_j}\int_{t_j}^ty dt \right|\leq \frac{\epsilon}{t-t_j}
\end{equation}

Consider $\frac{1}{t-t_j}\int_{t_j}^t x \ dt$ and note, as was discussed before, that 
\begin{align} 
x(t)-x(P)<\nu \ \Rightarrow x(P)\leq x(t)\leq x(P)+\nu
\end{align}
which yields  
\begin{align} 
x(P)\leq \frac{1}{t-t_j}\int_{t_j}^tx dt\leq x(P)+\nu.
\end{align}
We abbreviate $\bar{x}(j,t)=\frac{1}{t-t_j}\int_{t_j}^tx dt; \ \ \bar{y}(j,t)=\frac{1}{t-t_j}\int_{t_j}^ty dt$. Hence we find that $(\bar{x}(j,t),\bar{y}(j,t))$ lies in the $\nu$-neighbo  rhood of $P$.

Now choosing $\epsilon$ small enough and $t$ big enough, we find that the right hand side  $\frac{\epsilon}{t-t_j}$ of the inequality \eqref{eq:inq1} becomes arbitrarily small. This in turn implies that $\bar{x}(j,t), \ \bar{y}(j,t)$ closely satisfies the linear equation 
\begin{align} 
\frac{\alpha}{\beta}(1-\bar{x}(j,t))=\bar{y}(j,t),
\end{align}which is the line $g_x$. The analogous argument for $\frac{\dot{y}}{y}$ implies that $\bar{x}(j,t), \ \bar{y}(j,t)$ closely satisfies the equation of the line $g_y$, hence $\bar{x}(j,t), \ \bar{y}(j,t)$ are arbitrarily close to the critical point $C$. But we have shown that $\bar{x}(j,t), \ \bar{y}(j,t)$ converge to the point $P$, which is different from $C$ by assumption, hence for $\epsilon$ small enough and $t$ big enough we get the desired contradiction. 
\end{proof}

With these prerequisites we can now turn towards the classification:

\begin{theorem}\label{limit}
Every trajectory in the \emph{PC system} converges to one of the critical points.
\end{theorem}

\begin{proof}
Any trajectory of the \emph{PC system} is bounded since it is confined to $\Delta$, theorem \ref{ass} implies that the limit set of any point is then non-empty, compact and connected. Together with the fact that there exist only four critical points we can apply theorem \ref{PB}. Since theorem \ref{main} prohibits the existence of periodic orbits, the limit-set of any point in $\Delta$ is either a critical point or a connected set consisting of some critical points together with trajectories between them. But proposition \ref{lim} also excludes non-critical points on a trajectory connecting the critical points. Thus we are left with only critical points as limit-sets. The same argument as in the proof of proposition \ref{lim} guarantees that the convergence is actually true in the continuous sense (to exclude the possibility that a trajectory might recede from a limit point between times $t_i, t_j$). We briefly recall the argument. Either $\upsilon$ is completely contained in sectors $I$ or $III$, then we are either above or below both of $g_x,g_y$ and the convergence is monotonic. If it enters $II$ or $IV$, by proposition \ref{fish-trap} it will stay there for all times and also experience monotonic convergence. 
\end{proof}

We summarize what we have obtained:
Given any Point $Q\in \Delta$, if we follow the trajectory through $Q$ in the \emph{PC system} long enough, we will end up arbitrarily close to one of the critical points $C_0,C_x,C_y,C$. Which one is determined by the parameters $\alpha, \beta,\gamma,\delta$. Hence the long-term behaviour of the language change modelled by the \emph{PC system} is completely determined by the parameters. 

We are going to describe this dependence more precisely. To this end, consider the flow of the dynamical system in the vicinity of a critical point $C_p$ at a specific time $t_0$ and at a time $t_0+dt$ infinitesimally later. Since the flow fixes the critical point $C_p$ (the flow is stationary there), points near $C_p$ are flowing to points near $C_p$. Hence the flow between time $t_0$ and time $t_0+dt$ amounts to a linear map of neighborhoods of $C_p$, the \emph{linearisation} of the flow which captures its essential features. The linear map is given by the $2\times 2$ matrix $A$ of partial derivatives ("linearisation of the flow") of the equations of motion. Now a positive eigenvalue of $A$ corresponds to a point keeping its direction but flowing \emph{away} from the critical point (trajectories close to this one are thus repelled by the critical point), whereas a negative eigenvalue of $A$ corresponds to point keeping its direction but flowing \emph{towards} the critical point (thus nearby trajectories are attracted). To briefly describe the remaining cases: if an eigenvalue is zero, then the flow stagnates in the corresponding direction (consists of fixed points), whereas for a complex eigenvalues of $A$ the flow in the vicinity of the corresponding critical point would show rotatory character. Thus by calculating the eigenvalues of $A$ for all critical points of the \emph{PC system}, we can determine the behaviour of the flow.  
More details can be found in any standard textbook on dynamical systems or differential equations, for example in Zill and Cullen \cite{ZC}.  

Therefore we differentiate both defining equations of the \emph{PC system} with respect to $x,y$ and form the matrix:
\begin{align}
A=\left( 
\begin{array}{ll}
\dd{\dot{x}}{x} & \dd{\dot{x}}{y}\\ 
\dd{\dot{y}}{x} & \dd{\dot{y}}{y}
                             \end{array}\right)=\left(\begin{array}{ll}
\alpha -2\alpha x -\beta  y & -\beta x\\ 
-\delta y & \gamma -2\gamma  y-\delta  x 
                             \end{array}\right)
\end{align}
We need to substitute the critical points for $x,y$ into $A$ and then find the eigenvalues. A positive eigenvalue corresponds to a repelling eigendirection, a negative eigenvalue to an attracting eigendirection. 

\begin{enumerate}
 \item $C_0=(0,0)$: 
 \begin{align}
A=\left( \begin{array}{ll}
\alpha  & 0\\ 
0 & \gamma
\end{array}\right),
 \end{align}
obviously has two repelling eigendirections (the $x$- and $y$-axes).
\item $C_x=(1,0)$:
\begin{align}A=\left(\begin{array}{ll}
-\alpha  & -\beta\\ 
0 & \gamma-\delta
\end{array}\right),
\end{align}
the eigenvalues are $-\alpha$ and $\gamma-\delta$. 
By the discussion above $\gamma-\delta<0$. 
Thus there exist two attractive eigendirections (one along the $x$-axis). The critical point is thus a \emph{sink} and attracts all trajectories in the vicinity.
\item $C_y=(0,1)$:
\begin{align} 
A=\left(\begin{array}{ll}
\alpha-\beta  & 0\\ 
-\delta & -\gamma
\end{array}\right),
\end{align}the eigenvalues are $-\gamma$ and $\alpha-\beta$. This corresponds to two attractive eigendirections (one along the $y$-axis). The critical point is again a \emph{sink}. 
\item $C=(x_{crit},y_{crit})=(\frac{\gamma(\alpha-\beta)}{D},\frac{\alpha(\gamma-\delta)}{D})$:
\begin{align}
A=\left( \begin{array}{ll}
-\alpha \cdot x_{crit}  & -\beta\cdot x_{crit}\\ 
-\delta \cdot y_{crit} & -\gamma \cdot y_{crit}
\end{array}\right),
\end{align}

We begin by showing that the matrix $A$ must have a negative and a positive eigenvalue. If a $2\times 2$ matrix with real coefficients has a complex eigenvalue $\lambda$, it must have a second complex eigenvalue which equals the complex conjugate $\overline{\lambda}$ of the first, but then $det A=\lambda \cdot \overline{\lambda}>0$.
Now the determinant equals (cf. discussion on parameters above): 
\begin{align}
det A=D\cdot x_{crit}\cdot y_{crit}=\frac{\gamma(\alpha-\beta)\alpha(\gamma-\delta)}{D}<0.
\end{align}
Hence we must have a positive and a negative real eigenvalue of $A$. 
Therefore we must have a repelling and an attracting eigendirection. 
\end{enumerate}

\subsubsection{Description of the resulting scenario}

Please note that $C_0$ is a repelling critical point. Thus we discard it from the list below, because the trajectories will never end there. Further note that away from the $x,y$-axis, the trajectories are pushed into $\Delta$.
From the previous deduction, we find the following scenario: each of the critical points $C_x,C_y$ is attracting. Moreover $C$ is a saddle point (has a repelling and an attracting eigendirection). Thus except for two trajectories coming in exactly in the attracting eigendirection, all other trajectories end in either $C_x$ or $C_y$. These two trajectories constitute the so-called \emph{separatrix} of the \emph{PC system} since they separate the space of trajectories in those converging to $C_x$ or $C_y$.    

\subsection{Classification of long-term behaviour in the singular cases}\label{sing}
In this section we are going to repeat the analysis above for special cases of the interaction parameters and discuss their implications. Therefore recall the definitions
$$\alpha=N\widetilde{\alpha};$$ $$\beta=N (\widetilde{\alpha}+\widetilde{\beta}); \ \ \ \ \ \gamma=N \widetilde{\gamma}; \ \ \ \ \ \delta=N(\widetilde{\gamma}+\widetilde{\delta}).$$

In general, the critical point $C$ is situated in the interior $\mathring{\Delta}$ of the triangle $\Delta$. 
We call this the \emph{generic} situation. The situation where the critical point $C$ does not exist or is situated on the boundary $\partial \Delta$ of the triangle $\Delta$ is called \emph{singular}. The above analysis deals therefore with the \emph{generic} case. To analyse the \emph{singular} cases, we discuss the different possible loci which the critical point $C$ can have on the boundary $\partial \Delta$. Therefore consider Figure \ref{fig:3} for nomenclature:
\begin{figure}[ht!]
\psfrag{a}{$a$}
\psfrag{b}{$b$}
\psfrag{z}{$c$}
\psfrag{cy}{$C_y$}
\psfrag{cx}{$C_x$}
\psfrag{c}{$C$}
\psfrag{c0}{$C_0$}
 \includegraphics[width=8cm]{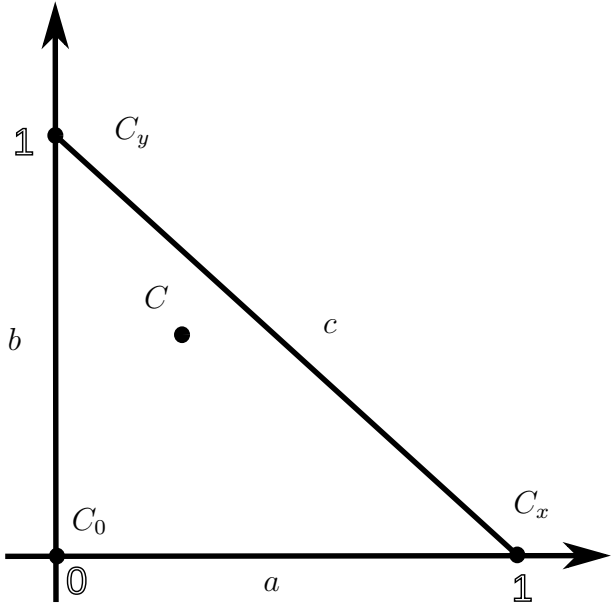}
 \caption{Nomenclature of $\partial \Delta$}
 \label{fig:3}
 
\end{figure}

\begin{enumerate}
 \item \emph{Generic case} \newline
 We included this point only to show that $C \in \mathring{\Delta}$ (our definition of being \emph{generic}) is equivalent to the conditions on the parameters used in the previous section. Indeed, from $C \in \mathring{\Delta}$ in follows straight away that: 
 \begin{align}
0<\frac{\gamma(\alpha-\beta)}{D}<1 \ \ \wedge \ \ 0<\frac{\alpha(\gamma-\delta)}{D}<1.
 \end{align}
 From this we directly conclude: $0<\gamma$ and $\alpha<\beta$ (by definition of $\alpha,\beta$ it is clear that $\alpha \leq \beta$) and analogously $0<\alpha$ and $\gamma<\delta$. Thus the generic case is equivalent to $0<\alpha<\beta$ and $0<\gamma<\delta$ together with $D<0$.
 \item \emph{One negative parameter}\newline
 It must be assured that the flow of the dynamical system points not outwards of $\Delta$. Otherwise the flow could generate a negative proportion of liberal speakers. The condition $\beta+\delta\geq \alpha+\gamma$ (by calculation of the dot-product $(\dot{x},\dot{y})^T \circ (-1,-1)^T\geq 0$) ensures this. A short inspection reveals that the critical point $C$ has at least one negative component, so that it is not situated within $\Delta$. Thus the flow of the dynamical system in $\Delta$ is reduced to two of the sectors $I,II,III,IV$ in Figure \ref{fig:1}. The conclusions are the same as before in such a situation (convergence to one of the critical points as in proposition \ref{lim}).
 \item \emph{$C=(q,0)\in a$}\newline
 The line $g_y: \ y=1-\frac{\delta}{\gamma}x$ is bound to go through $C_y=(0,1)$ whereas the line $g_x:\ y=\frac{\alpha}{\beta}(1-x)$ is bound to go through $C_x=(1,0)$. Since $C=g_x\cap g_y$, we see that $C=(q,0)$ forces $g_x$ to be identical with the x-axis. Hence $\alpha=0$ and $\dot{x}\leq 0$ in all of $\Delta$ (in our terminology from above, only the sectors $I, II$ survive in $\Delta$ since $III,IV$ lie beneath the line $g_x$). But on the x-axis ($y=0 \ \Rightarrow \ \ \dot{y}=0$) the flow can in general only be horizontal, but since $g_x$ is defined as the line where $\dot{x}=0$ vanishes, it follows straight away that all points on $a$ are critical. The methods of the generic case go through unaltered, except that we have to show that any limit-set can contain only one critical point (which follows straight away from a monotonicity argument). Thus we get the same conclusions:
 any trajectory ends at one of the critical points (here $C_y$ or a point of the segment $\{(p,0)| 0\leq p \leq 1\}$). As before consider the matrix $A$ obtained from linearisation of the flow. At the point $C_p=(p,0)$ we find: 
 \begin{align}
 A_p=\left( \begin{array}{ll}
0  & -\beta p\\ 
0 & \gamma- \delta  p
\end{array}\right).
 \end{align}
For $p<q$, the quantity $\gamma- \delta \cdot p$ is positive, whereas for $p>q$ it is negative. Hence to the left of the critical point $C$, the critical points of the $x$-axis are repelling whereas to the right of $C$ they are attracting. Again the fish-trapping lemma remains valid, thus any trajectory entering sector II can only converge to the critical point $C_y$, thus realising reversible language change as no progressive speakers exist. All other trajectories which never enter sector II will end up at one of the attractive critical points with $p>q$. This is an instance of incomplete language change since in the long term, we end up with a percentage $p$ of progressive speakers and the rest being liberal speakers (no conservative speakers left) not using the new feature.

\item \emph{$C=(0,q)\in b$}\newline
Exactly analogue to the previous case by exchanging $x,y$. The interpretation is as follows: any trajectory either ends on $C_x$, in the case of which we have complete language change, or they end at one of the critical points $(0,p)$ with $p>q$. In this case no progressive speakers are present and hence we have reversible language change (the new feature becomes extinct). 

\item \emph{$q \to 0$ in the third case}\newline
This means that $C \to C_0$ and that $\gamma \to 0$ (we already established $\alpha=0$ before). Thus we get the implications of both previous singular cases: All points on either $a,b$ are critical and attracting. A short calculation shows that the convergence to these critical points flows along the straight lines with equations:\\
$y=\frac{\delta}{\beta}\cdot x+ constant$. The interpretations from above carry over: either incomplete or reversible language change.
\item \emph{$q \to 1$ in fourth case}\newline
Now all trajectories are lying in sector $II$ and are thus converging to $C_y$. This is another instance of reversible language change with the new feature becoming extinct.
\item \emph{$q \to 1$ in third case}\newline
Now all trajectories are converging to $C_x$ and we find complete language change in all cases.
\item \emph{Points on $c$ become critical}\newline
In this case, we must have that both $g_x,g_y$ equal the line with equation $y=1-x$. This implies that $\alpha=\beta$ and $\gamma=\delta$. But then, all points on $y=1-x$ are critical (being on $g_x$ they satisfy $\dot{x}=0$ and by the same argument $\dot{y}=0$.) An argument similar to the one before shows that the points $C_p=(p,1-p)$ are attractive critical points and all trajectories will end at one of these critical points. Again this is an instance of incomplete language change, so that in this case we find only progressive (percentage $p$) and conservative speakers (percentage $1-p$) whereas no liberal speakers are present anymore.  
\end{enumerate}

\section{Application to general language change}
\subsection{General conclusions from the mathematical section}
In this section we only discuss the generic case, since for convenience we included the interpretation of the singular cases in section \ref{sing}.
The application to language change is now clear if we recall the meaning of the variable $x$ (the percentage of progressive speakers in the speaking community) and the variable $y$ (the percentage of conservative speakers in the speaking community).
Hence we find a very specific choice of initial data (ratios of progressive to conservative speakers), which results in an unstable equilibrium. Even the slightest alteration of this ratio leads to a totally different long-term behaviour, either to complete language change or to the extinction of the new feature (reversible language change).     
Thus we have now qualitatively solved the \emph{PLC model} for language change completely, given the interaction parameters $\alpha, \beta,\gamma,\delta$. Provided that the model proves successful to describe real data of language change, the gain from such a model is a considerable improvement to older approaches due to the fact that as soon as any data-fitting procedure has  
produced estimates for the interaction parameters, the \emph{PLC model} allows for predictions about the long-term behaviour of the system without having to do any further calculations. In the prevailing majority of modelled situations, fairly strong conclusions about the stability of the predictions can be concluded from the interaction parameters. 
\subsection{Testing the \emph{PLC model} on empirical data}\label{testing}
The aim of this section is to describe how we approached the numerical testing of the \emph{PLC model}.  
We are not going to repeat well-known facts about data-fitting procedures which have been described extensively, for example in the context of language change in Altmann \cite{A}. We use the relevant packages of the \emph{scipy, numpy} modules of the programming language \emph{Python}. To get better accessibility for the reader, we use a \emph{jupyter-notebook} which allows for a hybrid environment consisting of the executable Code and explanations thereof. The files used below can be downloaded from GitLab under the following link \url {https://gitlab.fosbos-rosenheim.de/pub/}.

Only one aspect which to our knowledge is not (obviously) standard material on data-fitting is the fact that due to the (presumed) non-integrability of the \emph{PLC model} we had to obtain the family of functions allowed for the data-fitting procedure by numerical integration. To this end we applied another Python module using a suitable \emph{Runge-Kutta solver}. We start the numerical analysis by comparing the \emph{PLC model} to the well-studied Piotrowski-Altmann law.

\subsubsection{Comparison to the well-studied Piotrowski-Altmann law}
In the following we compare the Piotrowski-Altmann law to the \emph{PLC model}. We take data from Best and Kohlhase on page 97 in \cite{BK}, which was also used by Altmann in \cite{A}. The data below shows the development over time of the percentage of usage of the (new) word \emph{wurde} opposed to the (old) word \emph{ward} in German.\\
\noindent\begin{minipage}{8.5cm}
\centering
\includegraphics[width=8.5cm]{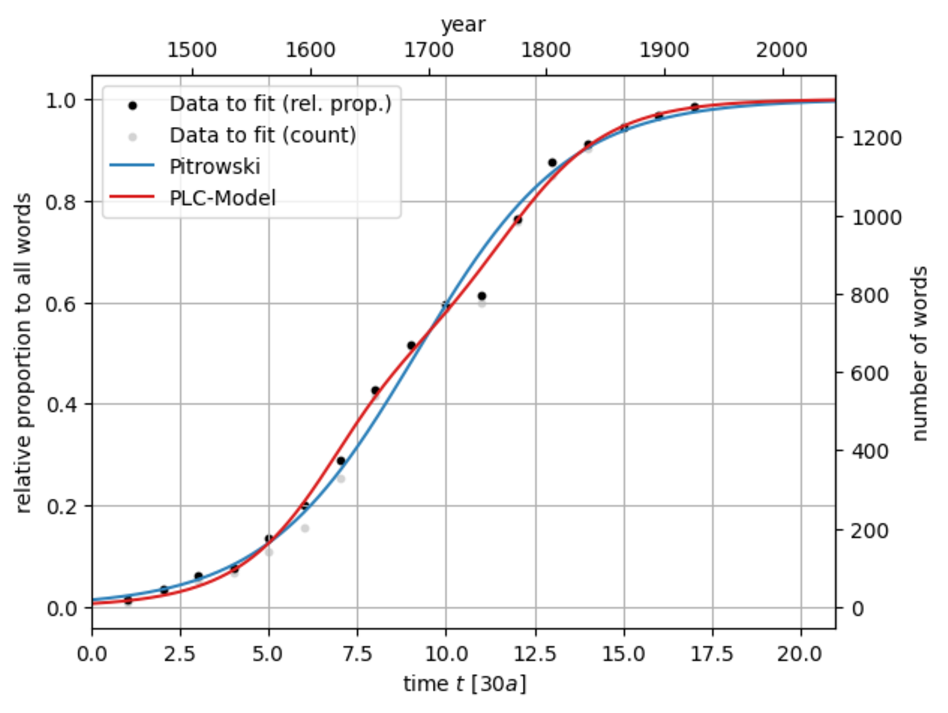}
\captionof{figure}{PLC and Piotrowski-Altmann}
\label{fig:1a}    
\end{minipage}\begin{minipage}{6cm}
Using Piotrowski's law, Altmann obtained as optimal function:
$f(t)=\frac{1}{1+70.4286\cdot e^{-0.4642\cdot t}}$. The graphic shows the solution from the \emph{PLC model} (red) in comparison to Altmann's optimal function (blue) on the dataset described above. The fitting parameters for both models can be found in the appendix \ref{Papio}.  
\end{minipage}

\subsubsection{Testing on the e-epithesis}
The \emph{e-epithesis} is, according to Imsiepen \cite{I}, an example of reversible language change. By e-epithesis one understands the phenomenon in Early Modern High German to put an additional \emph{-e} to the end of strong verbs in the preterite, for example \emph{sahe} opposed to \emph{sah}. Over time this tendency initially started to increase before eventually growing completely out of fashion. We use the data from Imsiepen \cite{I}.\\
\noindent\begin{minipage}{8.5cm}
\centering
\includegraphics[width=8.5 cm]{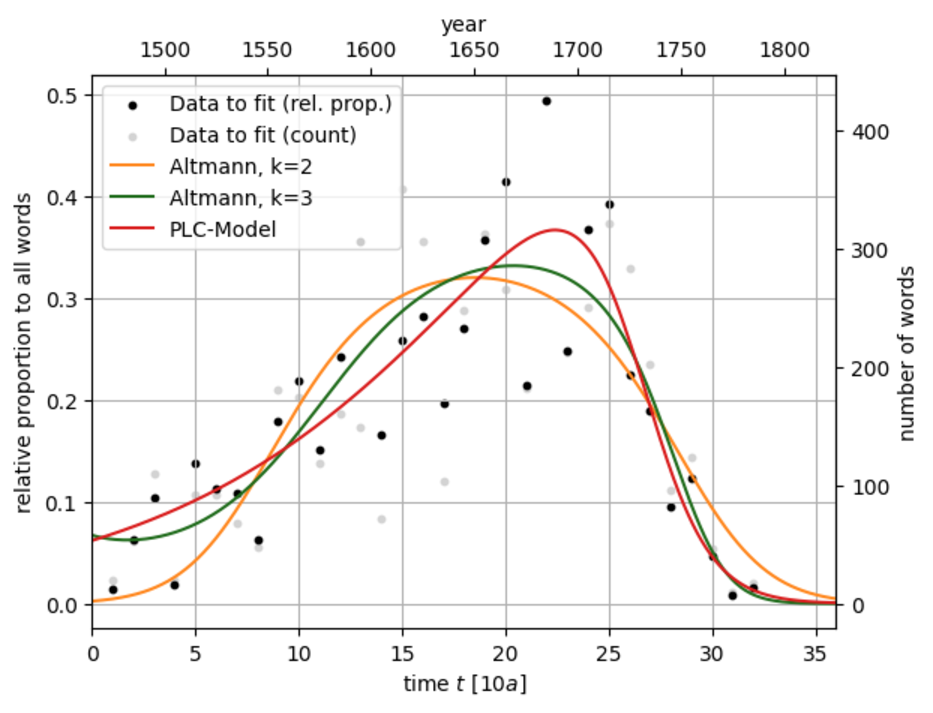}
\captionof{figure}{e-ephithesis}
\label{fig:5}
\end{minipage}\begin{minipage}{6cm}
In Figure \ref{fig:5}, we compare the fitting of the model from Altmann in \cite{A}, the generalized logistic curve from \cite{VB} using a third-order polynomial in the time variable and the PLC model. The graph of Altmann's attempt is orange, the generalized logistic curve is green whereas ours is red. The parameters of each model are given in the appendix \ref{Paeph}. Although the data is quite spread out, one can see that our model does show an improvement to the older models.
\end{minipage}

\subsubsection{Testing on the development of periphrastic \emph{do}}
According to Vulanovi\'{c} and Baayen \cite{VB}, there is another form of language change, exemplified by the proportion of periphrastic-\emph{do} constructions around 1560 which can be viewed as a two-stage change. By this it is meant that the data shows a slow-down (or even a decrease) before it starts increasing again. Vulanovi\'{c} and Baayen used a generalised logistic function in \cite{VB}, where they applied (besides another approach) a polynomial of order $k=3$ in the time-variable. We compare the generalized logistic function to the \emph{PLC model}. Even though the \emph{PLC model} cannot cope with an actual decrease (under the assumption of time-independent coefficients!), the following example shows how fitting of a slow-down is realised by the \emph{PLC model}. In \cite{VB} six different sentence types are analysed, here we consider only the two most interesting ones for our application and refer the interested reader to \url{https://gitlab.fosbos-rosenheim.de/pub/} for fitting of the other examples. Originally, the data has been collected by Elleg\aa{a}rd \cite{El}(Table 7). We use the data as presented in \cite{VB} (Table 1). There the  thirteen periods in \cite{El} have been reduced to eleven to compensate for the different number of texts in the periods considered (cf. \cite{VB} p. 2 for further discussion).  

\noindent\begin{minipage}{7.3cm}
\centering
\includegraphics[width=7.3cm]{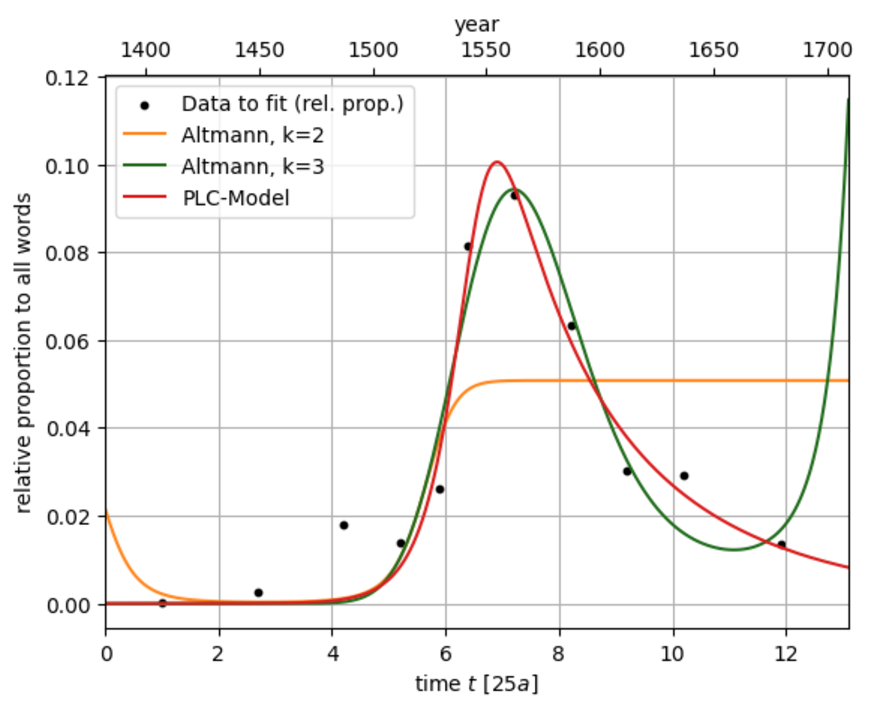}
\captionof{figure}{Affirmative declarative}\label{fig:7}
\end{minipage}\begin{minipage}{7.3cm}
\centering
\includegraphics[width=7.1cm]{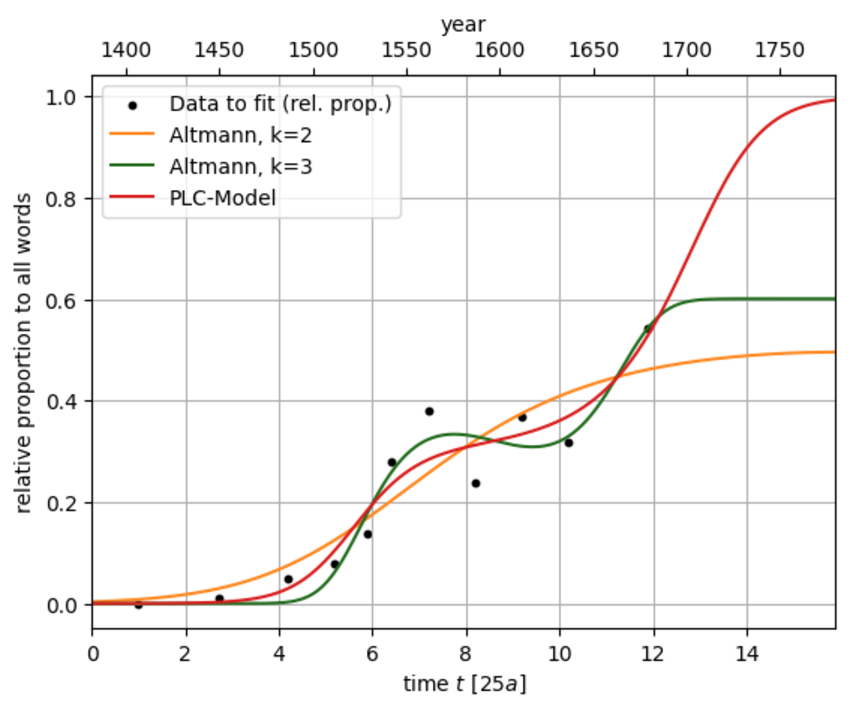}
\captionof{figure}{Negative declarative}\label{fig:9}
\end{minipage}

\vspace{0.5cm}
The graphs show good fitting results for both models (parameters are given in \ref{PaAD}), but the \emph{PLC model} shows less tendency to overshooting on the (time-)ends of the datasets (Figures \ref{fig:7}). 
\subsubsection{Prediction capabilities of the PLC model}
In this section we want to examine the capability of the \emph{PLC model} to predict the future development of processes of language change. We therefore left out the last six data points in the case of the e-ephithesis and the last three data points in the case of affirmative declarative \emph{do} and applied the fitting procedures to the remaining datasets (parameters are given in \ref{Paephpred}). In the graphics below the missed out data points are displayed in red so that it becomes possible to judge how well the model approximates the future development of the process under consideration.

\noindent\begin{minipage}{7.3cm}
\centering
\includegraphics[width=7.3cm]{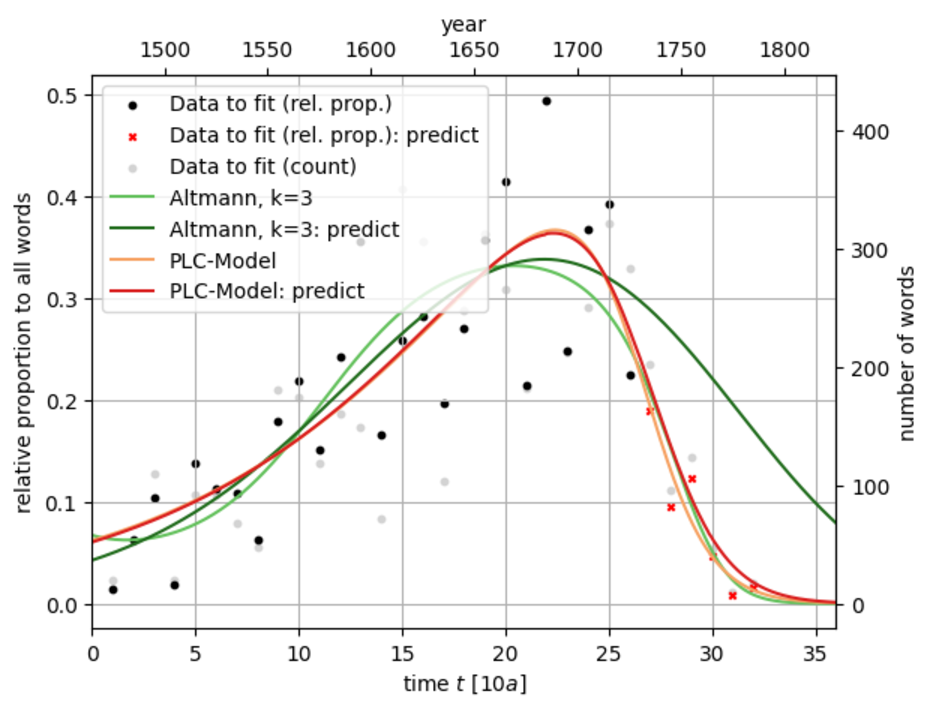}
\captionof{figure}{e-ephithesis}\label{fig:10}
\end{minipage}\begin{minipage}{7.3cm}
\centering
\includegraphics[width=6.9cm]{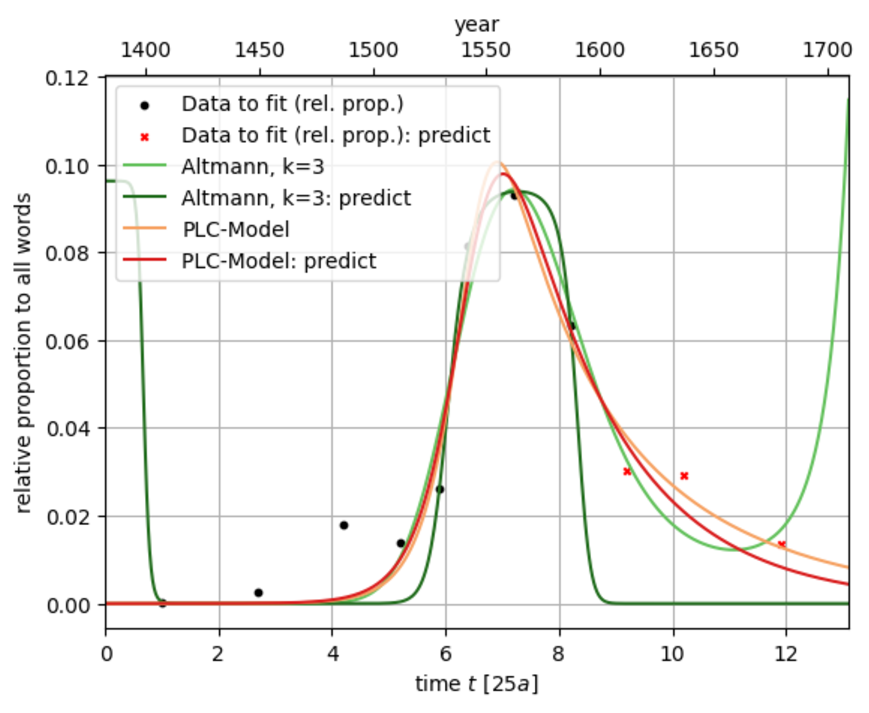}
\captionof{figure}{Affirmative declarative}\label{fig:11}
\end{minipage}
\section{Conclusion}
As stated, the \emph{PLC model} is not limited to the context of linguistics. Since it abstractly describes the interplay between \emph{progressive, liberal} and \emph{conservative} influences, it should be also applicable to various different settings situated in sociology, economics or political sciences. This is subject to further study. 

\section{Acknowledgement}
We are deeply grateful to Professor Gabriel Altmann, without whose positive feedback we would have never considered to prepare this work for publication. 
This article grew out of a semester work in 2020 of one of us (E.S.) and an email exchange with Prof. Altmann. Sadly, during the preparation, Prof. Altmann passed away.
We also like to thank the anonymous referees for valuable suggestions which considerably improved the content and the presentation of the manuscript.
\section{Appendix}
For the sake of completeness, we present the parameters of the fitting functions in section \ref{testing}. First recall the definition of parameters of the fitting functions:
\begin{enumerate}
\item \textbf{Piotrowski-Altmann}\\
$f(t)=\frac{c}{1+ae^{-bt}}$
 \item \textbf{Altmann, $k=2$}\\
 $f(t)=\frac{c}{1+ae^{-(bt+dt^2)}}$
 \item \textbf{Altmann, $k=3$}\\
 $f(t)=\frac{c}{1+ae^{-(k_0+k_1t+k_2t^2+k_3t^3)}}$
 \item \textbf{PLC model}\\
 $x(t)$ solving 
 \begin{align}
\dot{x}&=&\alpha  x (1-x)-\beta xy\\
\dot{y}&=&\gamma y(1-y)-\delta  x y
 \end{align}
\end{enumerate}
\begin{itemize}
 \item\textbf{Piotrowski-Altmann; Figure \ref{fig:1a} \label{Papio}}\\
 \begin{itemize}
 \item \textbf{Piotrowski-Altmann}\\
$\displaystyle{a = 70.4286 }$\\
$\displaystyle{b = 0.4642 }$\\
$\displaystyle{c = 1.0000 }$

\item \textbf{PLC model}\\
$\displaystyle{\alpha = 0.6142 \pm 0.1381}$\\
$\displaystyle{\beta = 2.6240 \pm 50.0777}$\\
$\displaystyle{\gamma = 2.1509 \pm 5.4645}$\\
$\displaystyle{\delta = 4.2129 \pm 7.2251}$\\
$\displaystyle{x_0 = 0.0067 \pm 0.0044}$\\
$\displaystyle{y_0 = 0.0000 \pm 0.0000}$
\end{itemize} 
 
 \item \textbf{e-epithesis; Figure \ref{fig:5} \label{Paeph}}\\
 \begin{itemize}
 \item \textbf{Altmann $k=2$}\\
$\displaystyle{a = 174.7431 \pm 323.4998}$\\
$\displaystyle{b = 0.6895 \pm 0.3641}$\\
$\displaystyle{c = 0.4139 \pm 0.1840}$\\
$\displaystyle{d = -0.0186 \pm 0.0098}$
\item \textbf{Altmann $k=3$}\\
$\displaystyle{c = 0.3973 \pm 0.1320}$\\
$\displaystyle{k_0 = -1.5936 \pm 1.2400}$\\
$\displaystyle{k_1 = -0.1060 \pm 0.3113}$\\
$\displaystyle{k_2 = 0.0336 \pm 0.0303}$\\
$\displaystyle{k_3 = -0.0010 \pm 0.0007}$
\item \textbf{PLC model}\\
$\displaystyle{\alpha = 0.1076 \pm 0.0256}$\\
$\displaystyle{\beta = 2.3732 \pm 144.0511}$\\
$\displaystyle{\gamma = 0.0377 \pm 1.4282}$\\
$\displaystyle{\delta = -1.1806 \pm 4.1077}$\\
$\displaystyle{x_0 = 0.0618 \pm 0.0199}$\\
$\displaystyle{y_0 = 0.0001 \pm 0.0034}$
\end{itemize}
\item \textbf{affirmative declarative (AD); Figure \ref{fig:7}\label{PaAD}}\\
 \begin{itemize}
 \item \textbf{Altmann $k=2$}\\
$\displaystyle{a = 1.3349 \pm 29.0541}$\\
$\displaystyle{b = -3.4284 \pm 11.6307}$\\
$\displaystyle{c = 0.0508 \pm 0.0123}$\\
$\displaystyle{d = 0.6198 \pm 1.5044}$
\item \textbf{Altmann $k=3$}\\
$\displaystyle{c = 4288.0894 \pm 281522975.0661}$\\
$\displaystyle{k_0 = -58.0758 \pm 65802.4601}$\\
$\displaystyle{k_1 = 16.7839 \pm 15.5784}$\\
$\displaystyle{k_2 = -1.9223 \pm 1.8064}$\\
$\displaystyle{k_3 = 0.0701 \pm 0.0671}$
\item \textbf{PLC model}\\
$\displaystyle{\alpha = 2.1602 \pm 2.3368}$\\
$\displaystyle{\beta = 2.5054 \pm 2.1382}$\\
$\displaystyle{\gamma = 3.5855 \pm 8.9612}$\\
$\displaystyle{\delta = -0.3424 \pm 12.2549}$\\
$\displaystyle{x_0 = 0.0000 \pm 0.0000}$\\
$\displaystyle{y_0 = 0.0000 \pm 0.0000}$
\end{itemize}
\item \textbf{negative declarative (ND); Figure \ref{fig:9}\label{PaND}}\\
 \begin{itemize}
 \item \textbf{Altmann $k=2$}\\
$\displaystyle{a = 159.0507 \pm 2610.5674}$\\
$\displaystyle{b = 0.8750 \pm 3.9980}$\\
$\displaystyle{c = 0.5660 \pm 4.0071}$\\
$\displaystyle{d = -0.0272 \pm 0.4874}$
\item \textbf{Altmann $k=3$}\\
$\displaystyle{c = 0.6009 \pm 0.2411}$\\
$\displaystyle{k_0 = -44.1164 \pm 33.5010}$\\
$\displaystyle{k_1 = 15.7658 \pm 13.1213}$\\
$\displaystyle{k_2 = -1.8539 \pm 1.6519}$\\
$\displaystyle{k_3 = 0.0720 \pm 0.0683}$
\item \textbf{PLC model}\\
$\displaystyle{\alpha = 1.7308 \pm 10.8865}$\\
$\displaystyle{\beta = 8.2517 \pm 366.9127}$\\
$\displaystyle{\gamma = 1.1318 \pm 17.8137}$\\
$\displaystyle{\delta = 3.1689 \pm 31.6397}$\\
$\displaystyle{x_0 = 0.0000 \pm 0.0010}$\\
$\displaystyle{y_0 = 0.0003 \pm 0.0160}$
\end{itemize}
\item \textbf{e-ephithesis-prediction; Figure \ref{fig:10}\label{Paephpred}}\\
 \begin{itemize}
 \item \textbf{Altmann $k=3$}\\
$\displaystyle{c = 0.3973 \pm 0.1320}$\\
$\displaystyle{k_0 = -1.5937 \pm 1.2400}$\\
$\displaystyle{k_1 = -0.1059 \pm 0.3113}$\\
$\displaystyle{k_2 = 0.0336 \pm 0.0303}$\\
$\displaystyle{k_3 = -0.0010 \pm 0.0007}$
\item \textbf{Altmann $k=3$-predict}\\
$\displaystyle{c = 1.9536 \pm 113.1033}$\\
$\displaystyle{k_0 = -3.8067 \pm 60.6009}$\\
$\displaystyle{k_1 = 0.1604 \pm 0.4607}$\\
$\displaystyle{k_2 = -0.0006 \pm 0.0969}$\\
$\displaystyle{k_3 = -0.0001 \pm 0.0026}$
\item \textbf{PLC model}\\
$\displaystyle{\alpha = 0.1076 \pm 0.0256}$\\
$\displaystyle{\beta = 2.3732 \pm 144.0511}$\\
$\displaystyle{\gamma = 0.0377 \pm 1.4282}$\\
$\displaystyle{\delta = -1.1806 \pm 4.1077}$\\
$\displaystyle{x_0 = 0.0618 \pm 0.0199}$\\
$\displaystyle{y_0 = 0.0001 \pm 0.0034}$
\item \textbf{PLC model-predict}\\
$\displaystyle{\alpha = 0.1109 \pm 0.0814}$\\
$\displaystyle{\beta = 2.7888 \pm 585.9505}$\\
$\displaystyle{\gamma = 0.0614 \pm 15.6047}$\\
$\displaystyle{\delta = -0.9330 \pm 43.0411}$\\
$\displaystyle{x_0 = 0.0602 \pm 0.0244}$\\
$\displaystyle{y_0 = 0.0001 \pm 0.0327}$
\end{itemize}
\item \textbf{affirmative declarative-prediction; Figure \ref{fig:11}\label{PaADpred}}\\
 \begin{itemize}
 \item \textbf{Altmann $k=3$}\\
$\displaystyle{c = 4288.0894 \pm 281522975.0661}$\\
$\displaystyle{k_0 = -58.0758 \pm 65802.4601}$\\
$\displaystyle{k_1 = 16.7839 \pm 15.5784}$\\
$\displaystyle{k_2 = -1.9223 \pm 1.8064}$\\
$\displaystyle{k_3 = 0.0701 \pm 0.0671}$
\item \textbf{Altmann $k=3$-predict}\\
$\displaystyle{c = 0.0962 \pm 0.0184}$\\
$\displaystyle{k_0 = 15.4560 \pm 50.7514}$\\
$\displaystyle{k_1 = -27.1872 \pm 30.2628}$\\
$\displaystyle{k_2 = 6.8091 \pm 6.0113}$\\
$\displaystyle{k_3 = -0.4529 \pm 0.3692}$
\item \textbf{PLC model}\\
$\displaystyle{\alpha = 2.1602 \pm 2.3368}$\\
$\displaystyle{\beta = 2.5054 \pm 2.1382}$\\
$\displaystyle{\gamma = 3.5855 \pm 8.9612}$\\
$\displaystyle{\delta = -0.3424 \pm 12.2549}$\\
$\displaystyle{x_0 = 0.0000 \pm 0.0000}$\\
$\displaystyle{y_0 = 0.0000 \pm 0.0000}$
\item \textbf{PLC model-predict}\\
$\displaystyle{\alpha = 2.0863 \pm 16.9947}$\\
$\displaystyle{\beta = 2.6160 \pm 197.3765}$\\
$\displaystyle{\gamma = 2.6927 \pm 222.9634}$\\
$\displaystyle{\delta = 1.1746 \pm 2518.9027}$\\
$\displaystyle{x_0 = 0.0000 \pm 0.0000}$\\
$\displaystyle{y_0 = 0.0000 \pm 0.0001}$
\end{itemize}

\end{itemize}

\end{document}